\documentclass[12 pt]{amsart}
\usepackage{amsfonts,amscd,amsmath,amssymb,amsthm}
\usepackage{times}
\usepackage[centertags]{amsmath}
\usepackage[dvips]{graphicx}

\DeclareSymbolFont{SY}{U}{psy}{m}{n}
\DeclareMathSymbol{\emptyset}{\mathord}{SY}{'306}

\DeclareMathOperator*{\slim}{s-lim} 
 
 \DeclareMathOperator{\tr}{tr}

 \DeclareMathOperator{\di}{dim}

\DeclareMathOperator{\les}{\lambda_{ess}}
\DeclareMathOperator{\diag}{diag}
\DeclareMathOperator{\spec}{\Sigma}

\DeclareMathSymbol{\newtimes}{\mathbin}{SY}{'264}

\newcommand{\fD}{\mathfrak{D}}

\newcommand{\fG}{\mathfrak{g}}

\newcommand{\cT}{\mathcal{T}}

\title{On Temple--Kato like inequalities and applications}
\author{Luka Grubi\v{s}i\'{c}}
\address{Institut f{\"u}r reine und angewandte Mathematik,
RWTH Aachen, Templergraben 55,
D-52056 Aachen, Germany.
(On leave from PMF-Department of Mathemathics, University of Zagreb, Croatia)}
\email{luka.grubisic@iram.rwth-aachen.de}
\thanks{This work is based on a part of author's PhD thesis \cite{GruPhd}, which was written
under the supervision of Prof. Dr. Kre\v{s}imir Veseli\'{c}, Hagen in partial
fulfilment of the requirements for the degree Dr. rer. nat.}

\begin{document}
%
%  defs
%
\def\ra{{\sf R}}
\def\je{{\sf N}}
\def\I{\mathbf{I}}
\def\mA{\mathbf{A}}
\def\mM{\mathbf{M}}
\def\mU{\mathbf{U}}
\def\mX{\mathbf{X}}
\def\mV{\mathbf{V}}
\def\mW{\mathbf{W}}
\def\mY{\mathbf{Y}}
\def\mK{\mathbf{K}}
\def\mR{\mathbf{R}}
\def\mQ{\mathbf{Q}}
\def\mT{\mathbf{T}}
\def\mS{\mathbf{S}}
\def\mH{\mathbf{H}}
\def\mB{\mathbf{B}}
\def\x{\mathcal{X}}
\def\q{\mathcal{Q}}
\def\b{\mathcal{B}}
\def\g{\mathcal{G}}
\def\d{\mathcal{D}}
\def\K{\mathcal{K}}
\def\H{\mathcal{H}}
\def\vp{\mathcal{V}}
\def\y{\mathcal{Y}}
\def\wp{\mathcal{W}}
\def\lp{\mathcal{L}}
\def\EE{\mathcal{E}}
\def\SS{\mathcal{S}}
\def\R{ \mathbb{R}}
\def\C{ \mathbb{C}}
\def\N{ \mathbb{N}}

%
% environments
%
%
% environments
%
\newtheorem{theorem}{Theorem}[section]
\newtheorem{corollary}[theorem]{Corollary}
\newtheorem{lemma}[theorem]{Lemma}
\newtheorem{proposition}[theorem]{Proposition}
\theoremstyle{definition}
\newtheorem{remark}[theorem]{Remark}
\newtheorem{assumption}[theorem]{{\rm \textbf{Assumption}}}
\newtheorem{example}[theorem]{Example}
\newtheorem{problem}{{\rm \textbf{Problem}}}[section]
\newtheorem{definition}[theorem]{Definition}
\numberwithin{equation}{section}
%
% commands
%
\newcommand{\norm}[1]{\Vert #1 \Vert}
\newcommand{\normm}[1]{\Vert #1 \Vert}
\newcommand{\ugl}[1]{\left[ #1 \right]}
\newcommand{\sk}[1]{\left( #1 \right)}
\newcommand{\dual}[1]{\left< #1 \right>}
\newcommand{\abs}[1]{\left| #1 \right|}
\newcommand{\abss}[1]{\vert #1 \vert}
\newcommand{\spa}[1]{\mbox{span}\{ #1 \}}
\newcommand{\absf}[2]{\frac{\abss{#1}}{\abss{#2}} }
\newcommand{\ine}[1]{{\mathbf{#1}}}
\newcommand{\conpo}[1]{\stackrel{#1}{\longrightarrow}}
%\definecolor{myblue}{rgb}{0,0,0.8}
%
% functions
%
\def\imag{{\rm i}}
\def\sinbf{\mathsf{sin}}
\def\cosbf{\mathsf{cos}}
\def\eexp{\text{e}}
\def\region{\mathcal{R}}
\def\slim{\text{{\rm s-lim}}}
\def\wlim{\text{{\rm w-lim}}}
\def\tripleb{ \mid\!\mid\!\mid }
%{\color{myblue}

\begin{abstract}
We give both lower and upper estimates for eigenvalues of unbounded
positive definite operators in an arbitrary Hilbert space. We show
scaling robust relative eigenvalue estimates for these operators in
analogy to such estimates of current interest in Numerical Linear
Algebra. Only simple matrix theoretic tools like Schur complements
have been used. As prototypes for the strength of our method we
discuss a singularly perturbed Schroedinger operator and study
convergence estimates for finite element approximations. The
estimates can be viewed as a natural quadratic form version of the
celebrated Temple--Kato inequality.
\end{abstract}
%\ams{
%34L15, 35P15, 65F15, 49R50}
%\end{AMS}
\maketitle
%\pagestyle{myheadings}
%\thispagestyle{plain}
%\markboth{L. Grubi\v{s}i\'{c}, \today }{L. Grubi\v{s}i\'{c}, \today}
%\markboth{On Temple--Kato like inequalities and applications}%
%{L. Grubi\v{s}i\'{c}}

\section{Introduction}
The purpose of this article is to establish scaling robust estimates
for discrete eigenvalues of positive definite operators in a Hilbert space.
We also prove that our estimates are optimal for
a residual type analysis of lowermost eigenvalues of those operators. Our approach uses
the theory of quadratic forms from \cite[Chapters VI--VIII]{Kato76}
and an adaptation of the matrix relative perturbation theory. As a
result we establish the same high performance residual type
estimates from \cite{DrmHAri97} in our more general setting. For a
review of the matrix relative perturbation theory see \cite{Li-I-98}
and the references therein.

It turns out that positive-definiteness of the matrix is the key
structural property which is needed for the analysis of
\cite{DrmHAri97,DrmacVeselic2}. Subsequently, we prove our estimates
for an abstract positive definite form in an arbitrary Hilbert
space. This, together with the fact that the estimates also hold for
the discrete eigenvalues which are in gaps of the essential
spectrum, indicates that our simple matrix analytic techniques are
well adapted to the class of problems under study, e.g. our
technique yields high performance estimates without forcing us to impose any unnecessary restrictions.
This abstract approach is further justified by the fact that we
simultaneously consider applications of these estimates to a study
of the convergence properties of adaptive finite element methods as
well as to a quantitative study of the asymptotic properties of
eigenvalue problems in the large coupling limit. Typical
operators in the large coupling limit setting are those from
\cite{BruneauCarbou,DemuthJeskeKirsch,SanchezPalencia90}.

The obtained estimates are the same as those which have proven
themselves in \cite{DrmacVeselic2} as a significant tool in the
development of modern mathematical software. The main feature of the
matrix eigenvalue algorithms from \cite{DrmacVeselic2} is that they
are robust when applied to extremely badly scaled input matrices. We
bring this in correspondence with the behavior of the spectrum of
stiffly/singularly perturbed operators from
\cite{BruneauCarbou,DemuthJeskeKirsch,SanchezPalencia90}. We use a
model Schroedinger operator to show that our estimates are optimal
for the class of stiff/singular perturbations, see Section
\ref{Section5}. An extensive study of non-inhibited stiff families
of operators has been performed, with the help of the results from
this article, in \cite{GruPhd} and will be published in the
subsequent report.

When studying the convergence of finite element procedures we adopt
the approach of \cite{Arioli,Neymeyr02}. Our results give a new
flavor to the analysis of the stoping criteria for preconditioned
inverse iterations from \cite[Section 4]{Neymeyr02}. In comparison,
we are more explicit about the dependence of the ingredients of the
error on the input data and we can prove the equivalence of the
error and the estimator, see Section \ref{Section6}. Furthermore, in
the Conclusion we briefly outline a simple way to obtain optimal
eigenvector estimates from our eigenvalue results.

In Section \ref{Section11} we relate our new results---which appear
in Section \ref{Section3}---to other results in the literature. In
Section \ref{Section12} we use a simple matrix model problem to give
a first flavor of the results by comparing our approach with that of
the Temple--Kato inequality.

\subsection{A comparison with other approaches}\label{Section11}
A number of recent studies of the eigenvalue approximation problem---when classified from the viewpoint of
perturbation theory (see e.g. \cite{Ipsen,Li-I-98})---could be seen to fall into the following
classes:
\begin{enumerate}
\item
Results obtained in the \textit{absolute}(sometimes also referred to as
regular) setting. For recent results see
\cite{Ovtchinik1} and references therein. The performance of such
estimates, when applied to an unbounded operator, depends on the
method of the regularization. The delicacy of this issue is
illustrated on an example below.
\item
Results which are obtained by interpreting the eigenvalue problem
as a nonlinear problem (in both eigenvalue and eigenvector). This
approach makes a treatment of the eigenvalue multiplicity somewhat
more difficult, since it is not easy to profit from the special
structure which the eigenvalue problem has, cf. \cite[Remark
7]{Heuveline}. %and \cite[Chapter 3.4]{Rudiger96}.
\item
Direct analysis of a representation of the single vector
\textit{residual} $r=\mH \psi-\mu \psi$, coupled with a consideration of the
approximation properties of the function space which is used to
generate the \textit{test vector} $\psi$, $\|\psi\|=1$, where we
have used $\mu:=(\psi, \mH \psi)$ to denote the \textit{Ritz value}.
Such estimates are essentially asymptotic in nature (cf.
\cite[Remark 3.2]{Duran}) or are specifically tailored for the
particular (class of) problem(s) under study (cf.
\cite[Theorem 4.1]{Larson2000}).
\end{enumerate}
We propose a technique which is based on the \textit{relative}
perturbation theory for qua\-dra\-tic forms in a Hilbert space from
\cite[Chapters VI--VIII]{Kato76}. In this regularization framework (which was developed in \cite{Gru03_3,GruVes02})
we solely use elementary matrix techniques, like LU-decomposition
and Schur complements, to obtain new block operator residual
equation. This residual equation has the same form as the
corresponding matrix result from \cite{DrmHAri97}, but holds in this
more general setting. It appears to be better suited to dealing with
eigenvalue multiplicity than are other approaches. We also argue
that when dealing with the lower part of the spectrum of the
positive definite operator our choice of the regularization is
optimal.

Specifically, we follow
 the approach of \cite{DemVes,RitzDrm96,DrmHAri97} and  reverse the trend to
  show that these finite dimensional results
have a lot to offer in the original setting of \cite{Kato76}, and in
particular as tools for a numerical study of singularly perturbed
(integro) differential operators. This paper will heavily use the
general construction from \cite{Gru03_3,GruVes02}. In short we
propose, assuming we are given a positive definite and self-adjoint
operator $\mH$ and the orthogonal projection $P$:
\begin{itemize}
\item
Construct the positive definite operator $\mH'$ from
\cite{Gru03_3,GruVes02}. $\mH'$ is the block diagonal part of $\mH$
with respect to $P$ and it is given formally as $\mH'=P\mH
P+(\I-P)\mH(\I-P)$. In this paper we shall use the notation
$\mH_P=\mH'$ to emphasize the dependence on the projection $P$.
\item
 Scale the perturbation $\mH-\mH_P$ with $\mH_P^{-1/2}$ to
obtain the bounded operator $\delta
H^P_s:=\text{``}\mH^{-1/2}_P(\mH-\mH_P)\mH^{-1/2}_P\text{''}$. This
is equivalent to working with the perturbation $\mH-\mH_P$ in the dual
energy space which is associated to $\mH_P$.
\item
Apply an adapted result from \cite{DemVes,RitzDrm96,DrmHAri97} to
obtain the desired eigenvalue/ vector estimate in the quadratic form setting.
\end{itemize}

At the end, we would like to emphasize that the original matrix
inequalities from \cite{DrmHAri97,DrmacVeselic2} have been extensively
tested in the process of developing new finite precision
eigenvalue software. In the course of this testing a large body of
test examples has been generated by judicious random searches as
well as by a modification of the known examples from science and
engineering. The inequalities have been found to be numerically sharp, as
is reported in \cite{DrmacVeselic2}, on numerous test matrices.

\subsection{Relationship to Temple--Kato inequality}\label{Section12} The central
theme of this paper is the issue of how to regularize an unbounded
eigenvalue problem to obtain an object which can then be
algebraically studied. As an introduction to the issue we shortly
review other approaches with special emphasis on the regularity
issues. This section is meant to provide the motivation for this
study and it extends the introduction.

The history of \textit{a posteriori} eigenvalue approximation
estimates goes back to
\cite{KatoT,Temple}. Such inequalities (most recently studied in
\cite{Ovtchinik1}) have a general form of
\begin{equation}\label{classical}
|\textsc{error}|\leq\textsc{conditioning}\times\|\textsc{residual}\|^2.
\end{equation}
Our estimates will have the form of
\begin{equation}\label{general_form}
\tripleb \textsc{residual}
\tripleb^2_{\textrm{rel}}\leq|\textsf{rel}(\textsc{error})|
\leq\textsf{rel}(\textsc{conditioning})\times\tripleb
\textsc{residual}\tripleb^2_{\textrm{rel}},
\end{equation}
where the measures $\tripleb.\tripleb_{\textrm{rel}}$ and
$|\textsf{rel}(\textsc{error})|$ and
$\textsf{rel}(\text{conditioning})$ denote the appropriate
ingredients from matrix (relative) perturbation theory as given in
\cite{Li-I-98}.

Let us now be more precise. We shall always work in the background
Hilbert space $\H$, which is
equipped with the scalar product $(\cdot,\cdot)$ and the norm
$\|\cdot\|=(\cdot,\cdot)^{1/2}$. Let $\mH$ be a self-adjoint
operator which is bounded from below and let $\psi$ be some vector
of the norm one in its \textit{domain of definition} $\d(\mH)$.
Define the \textit{Rayleigh-quotient} $\mu:=(\psi, \mH\psi)$ and
assume\footnote{We are counting the eigenvalues, which are below the infimum of the essential spectrum,
in the ascending order according to multiplicity.} that $\lambda_2(\mH)>\mu$ then
classical Temple--Kato inequality from \cite[Theorem VIII.5, Volume
IV pp. 84]{ReedSimonSvi} reads
\begin{equation}\label{TK:ineq}
\mu-\frac{(\mH\psi,\mH\psi)-(\psi,\mH\psi)^2}{\lambda_2(\mH)-\mu}\leq
\lambda_1(\mH)\leq\mu. \end{equation}

The vector $r=\mH\psi-\mu\psi$ is called the residual (associated to
$\psi$) and it holds $(\mH\psi,\mH\psi)-(\psi,\mH\psi)^2=\|r\|^2$.
Now it is easy to see that (\ref{TK:ineq}) implies
\begin{equation}\label{TK:formal}
\mu-\lambda_1\leq\frac{\|\mH\psi-\mu\psi\|^2}{\lambda_2-\mu},
\end{equation}
which has a general form of (\ref{classical}). Here we have used
$\lambda_i=\lambda_i(\mH)$, $i=1,2$ to simplify the notation.
The norm $\|r\|$ can be seen as an
``approximation defect'' of the vector $\psi$, since $\psi$ is an
eigenvector if and only if $\|r\|=0$.

As already stated our
estimates have a similar general form (see (\ref{general_form})) to
the Temple--Kato inequality but are obtained under the assumptions
of the perturbation theory for symmetric forms from \cite[Chapters
VI--VIII]{Kato76}. A consequence of this is that, in a case of a
positive definite operator $\mH$, we are able to directly work with test vectors
$\psi$ from the domain of the symmetric form which is according to
\cite[Theorem VI-2.23, pp. 331]{Kato76} equal to $\d(\mH^{1/2})$.
Our version of (\ref{TK:ineq}) also assumes $\mu<\lambda_2(\mH)$ but
allows $\psi\in\d(\mH^{1/2})$, $\|\psi\|=1$ and establishes (Theorem
\ref{thm:lower}) the estimate\footnote{The conditioning constant $\frac{\lambda_{2}+\mu}{\lambda_{2}-\mu}$
is a deliberate overestimate of the optimal constant $\fG_1$ from Theorem \ref{thm:lower}. It is a classical
ingredient of the \textit{relative} perturbation theory, see \cite{DrmacVeselic2,Li-I-98}.}
\begin{equation}\label{eq:ourTK}
\frac{\|\mH\frac{1}{\mu}\psi-\psi\|_{\mH^{-1}}^2}{\|\psi\|_{\mH^{-1}}^2}\leq\frac{\mu-\lambda_1}{\mu}\leq
\frac{\lambda_{2}+\mu}{\lambda_{2}-\mu} \frac{\|\mH\frac{1}{\mu}\psi-\psi\|_{\mH^{-1}}^2}{\|\psi\|_{\mH^{-1}}^2},
\end{equation}
where $\|\cdot\|_{\mH^{-1}}:=\|\mH^{-1/2}\cdot\|$ is the classical
$\mH^{-1}$-norm. To recognize the importance of the original
Temple--Kato approach, as well as in line with the terminology from
Numerical Linear Algebra, see \cite[pp. 271]{Neymeyr02}, we call all
the inequalities which have the form of (\ref{classical}) or
(\ref{general_form}) \textit{Temple--Kato like inequalities}. By
Temple--Kato approach we mean the notion that high performance
eigenvalue estimates should be obtained as a mixture of the
\textit{a posteriori} computable measure of the approximation defect
$\|r\|^2$ and the \textit{a priori} assumed quantitative information
$\displaystyle 1/(\lambda_2(\mH)-\mu)=\max
\{|\lambda-\mu|^{-1}~:~\lambda\in\spec(\mH)\setminus\{\lambda_1(\mH)\}\}$
on the conditioning of $\lambda_1(\mH)$.

Let us now discuss (\ref{eq:ourTK}). The measure of
$\textsf{rel}(\textsc{conditioning})$ is in this context the so called
\textit{relative gap} $(\lambda_{2}-\mu)/(\lambda_{2}+\mu)$, which
distinguishes close eigenvalues better than does the \textit{absolute
gap} $(\lambda_2-\mu)$ from (\ref{TK:ineq}). Furthermore, both the
residual measure $\|\mH \frac{1}{\mu}\psi-\psi\|^2_{\mH^{-1}}/\|\psi\|^2_{\mH^{-1}}$ as well as the
relative gap are robust with regard to scaling (e.g. ``dimensionless
quantities''). Thus, the most important message of (\ref{eq:ourTK}) is the
same as in \cite[Example 2.1]{DrmacVeselic2}: The
approximation $\mu$ has completely resolved the eigenvalue $\lambda_1$
when the (relative) residual measure drops below the relative gap. We
also note that (\ref{eq:ourTK}), unlike (\ref{TK:ineq}) cannot give
negative lower bounds to eigenvalues of positive definite
operators. On the other hand, the $\mH^{-1}$-norm is more difficult
to evaluate than are the ingredients of (\ref{eq:ourTK}). For a
possibility to do this see Remark \ref{rem:Golub}, Section
\ref{Section5} and \cite{BruneauCarbou,DemuthJeskeKirsch,GruPhd}.
Approximations to $\mH^{-1}$ norm of the residual can also be computed in a more accessible
scalar product, see \cite[Remark 7]{RitzDrm96}.
Note that the restriction $\psi\in\d(\mH)$ from (\ref{TK:ineq})
excludes---without prior regularization of the problem---the case when $\psi$ is a continuous piecewise linear
function and $\mH=-\triangle$ is the negative \textit{Laplace
operator}. On the other hand,  in this case is our theory directly
applicable and working estimates are explicitly given, accompanied with an argument for their optimality. Furthermore, as an illustration of our matrix theoretic
approach to unbounded operators, we will show (in Section
\ref{section2}) a ``matrix analytic'' way to obtain a variant of the
original Temple--Kato inequality.

We close this section by a simple and small numerical example which
should illustrate the dichotomy between the (easy) computability and scaling robustness of
eigenvalue estimates. We will be comparing the
first order estimates (in the approximation defect) from
\cite{Gru03_3,GruVes02} with the second order estimate (\ref{TK:formal}). As a model we consider the
asymptotic behavior of the family of positive-definite matrices
\begin{equation}\label{OvtEx}
\mH_\kappa=\begin{bmatrix} \frac{1}{101} & 0 & -\frac{1}{101}  \cr 0
&
    \frac{1}{100} & 0 \cr -\frac{1}{101}  & 0 & 1 +
   {\kappa }^2 \end{bmatrix},\qquad \kappa\to\infty.
\end{equation}

The results of \cite{Gru03_3,GruVes02}, which in this specialization
to the matrix $\mH_\kappa$ can be obtained using \cite[Theorem
1.1]{DrmHAri97} and a direct computation, use the relative residual
measure $$ \eta(\psi):=\frac{\|\mH\frac{1}{\mu}\psi-
\psi\|_{\mH^{-1}}}{\|\psi\|_{\mH^{-1}}},\qquad \mu=(\psi, \mH\psi).$$ Before we proceed note the
following geometrical facts. It holds that $\eta(\psi)$ is equal to
the sine of the angle $\angle(\mH \psi, \psi)$ in the scalar
product $(\cdot,\mH^{-1}\cdot)$. We denote this by
writing $ \eta(\psi)=\sin\angle(\mH \psi, \psi)_{\mH^{-1}} $. It can also
be shown that in the scalar product of the background space $\H$ the
identity $\eta(\psi)=\sin\angle(\mH^{1/2}\psi,\mH^{-1/2}\psi)$
holds. Furthermore, $\eta(\psi)$ is also an approximation defect
measure, since $\eta(\psi)\geq 0$ and $\psi$ is an eigenvector of
$\mH$ if and only if $\eta(\psi)=0$.

Now, \cite[Theorem 5.1]{Gru03_3} states that if for
$\psi\in\d(\mH^{1/2})$ the assumption
$\eta(\psi)<\frac{\lambda_2(\mH)-\mu}{\lambda_2(\mH)+\mu}$ holds
then
\begin{equation}\label{TK-gruves}
\frac{|\lambda_1(\mH)-\mu|}{\mu}\leq\eta(\psi).
\end{equation}
This estimate is equivalent with
\begin{equation}\label{TK-gruves2}
\big(1-\eta(\psi)\big)\mu\leq\lambda_1(\mH)\leq\big(1+\eta(\psi)\big)\mu.
\end{equation}
and since both $\mu$ and $\eta(\psi)$ are computable relations
(\ref{TK-gruves}) and (\ref{TK-gruves2}) give both a \textit{lower}
as well as an \textit{upper} estimate for $\lambda_1(\mH)$.

Now, take $\psi=\begin{bmatrix}1&0&0\end{bmatrix}^*$ as the test
vector and compare $\|r_\kappa\|$, $r_\kappa=\mH_\kappa\psi-\mu\psi$
and $\eta_\kappa(\psi)$,
$\eta_\kappa(\psi):=\|\mH_\kappa\frac{1}{\mu}\psi-
\psi\|_{\mH^{-1}_\kappa}/\|\psi\|_{\mH^{-1}_\kappa} $.

One computes $ \|r_\kappa\|=\frac{1}{101}$ whereas  $
\eta_\kappa(\psi)=\frac{1}{\kappa}\frac{\sqrt{2}}{\sqrt{101\kappa^{-2}+100}}$. This shows that the second order
estimate from (\ref{TK:formal}), as opposed to
the first order estimate (\ref{TK-gruves}), does not detect that
$(\mu-\lambda_1(\mH_\kappa))/{\mu}=\frac{1}{101
\kappa^2}+O\big(\frac{1}{\kappa^4}\big)\to 0$ as $\kappa\to\infty$. For more details on a numerical comparison
of (\ref{TK:ineq})--(\ref{TK:formal}) and (\ref{TK-gruves})--(\ref{TK-gruves2}) see \cite[Table 1.1]{Gru03_3}.

Caution has to be exercised when comparing absolute and relative estimates on this example.
It is known that absolute and relative estimation theory can (sometimes) yield
equivalent estimates, cf. \cite{Ipsen}. For instance, in the case of a single lowermost eigenvalue
$\lambda_1(\mH)$ an inequality which has a similar form as the righthand side inequality from
(\ref{eq:ourTK}) can be obtained if one applies (\ref{TK:formal}) to
the operator $\mH^{-1}$ in the Hilbert space with the $\mH$-scalar
product $(\cdot, \cdot)_\mH:=(\mH^{1/2}\cdot, \mH^{1/2}\cdot)$.
However, such approaches which first derive the eigenvalue estimates
in the background (absolute) scalar product and then scale the
operator at hand to fit this framework do not provide a proof of
the optimality of the estimates for parameter dependent
problems. Furthermore, our approach is more natural for the
treatment of the eigenvalue multiplicity which can be seen on the new
block-operator residual equation which yields error estimates that
utilize any unitary invariant norm of the block operator residual.
Such \textit{a posteriori} estimates and block operator residual equations did not appear before in the context of the
eigenvalue estimation for unbounded operators. We also note that the
numerical examples, reported in \cite{DrmacVeselic2}, indicate that
our choice of the (relative) residual measure $\eta(\psi)$, as well
as the choice of the measure of the conditioning (relative gap)
yield numerically sharp eigenvalue estimates. Bridging the relative
perturbation theory from \cite{DemVes,RitzDrm96,DrmHAri97} with the
theory of eigenvalue estimation for unbounded operators is the
declared aim of this work. In addition to that we will outline a
new general framework for analyzing asymptotic exactness of eigenvalue
estimators for parameter dependent eigenvalue problems. On our
simple $3\times 3$ example this general result reads
$$\lim_{\kappa\to\infty}\frac{\frac{\mu-\lambda_1(\mH_\kappa)}{\mu}}{\eta_\kappa^2(\psi)}=1$$
and the rate of the convergence appears to be rather rapid. We will
also show, by comparing the block operator residual equation which
yields (\ref{TK:formal}) with the residual equation that yields
(\ref{eq:ourTK}), that the approach of the relative perturbation
theory---e.g. first scale and then estimate rather than as in the
absolute approach where one first estimates and then scales---is the
right one when estimating the lower part of the
spectrum of a positive definite unbounded operator.

\section{A perturbation approach to Rayleigh--Ritz estimates}\label{section2}
We follow the general notational conventions and the terminology of
\cite[Chapters VI--VIII]{Kato76}. Minor differences are contained in the following list
of notation and terminology.
\begin{itemize}
\item $\H$ ... is an infinite dimensional Hilbert space, can be both real
  or complex
\item $(\cdot , \cdot)$; $\|\cdot\|$ ... the scalar product on $\H$, linear in the second
  argument and anti-linear (when $\H$ is complex) in the first; the norm on $\H$
\item $\H_1\oplus\H_2$... the direct sum of the Hilbert spaces $\H_1$ and $\H_2$, for any $x\in\H_1\oplus\H_2$
we have $x=x_1\oplus x_2=\begin{bmatrix}x_1\\x_2\end{bmatrix}$ for $x_i\in\H_i$, $i=1,2$
\item $\spec(\mH)$, $\spec_{ess}(\mH)$; $\les(\mH)$ ... the spectrum and the essential
  spectrum of $\mH$; the infimum of the essential spectrum of $\mH$
\item $A\leq B$ ... order relation between self-adjoint operators (matrices),
is equivalent with the statement that $B-A$ is positive
\item
$\lp(\H)$; $\lp(\H_1,\H_2)$... the space of bounded linear operators
on $\H$, which is equipped with the norm $\|\cdot\|$; the space of bounded linear
operators from $\H_1$ to $\H_2$
\item $\textsf{R}(X), \textsf{N}(X)$ ... the range and the null space of the linear operator $X$
\item $P$, $P_\perp$... the orthogonal projections $P$ and $P_\perp:=\I-P$
\item $j_{(\dot)}$ ... a permutation of $\N$
\item $\diag(M,W)$ ... the block diagonal operator matrix
with the operators $M,W$ on its diagonal. The operators $M,W$ can be
both bounded and unbounded. The same notation is used to define the
diagonal $m\times m$ matrix\\ $\diag(\alpha_1,\cdots,\alpha_m)$,
with $\alpha_1,\cdots,\alpha_m$ on its diagonal.
\item $s_1(A)\geq s_2(A)\geq\cdots$, $s_{\max(A)}, s_{\min}(A)$ ...
the singular values of the compact operator $A$ ordered in the
descending order according to multiplicity, the minimal (if it
exists) and the maximal singular value of $A$
\item $\tripleb X \tripleb$ ... a unitary invariant or operator cross norm
of the operator $X$. Since $\tripleb\cdot\tripleb$ depends
only on the singular values of the operator, we do not notationaly
distinguish between the instances of the norm
$\tripleb\cdot\tripleb$ on $\lp(\H)$, $\lp(\ra(P))$,
$\lp(\ra(P),\ra(P)^\perp)$, or such. Precise properties of a unitary
invariant norm will be listed in Section \ref{Section3}, for further details see \cite{SimonTrace}.
\item $\tr(X)$, $\tripleb X\tripleb_{HS}$ ... the trace a the Hilbert--Schmidt
norm of the operator $X$, it holds $\tripleb X\tripleb_{HS}=\sqrt{\tr(X^*X)}$, see \cite{SimonTrace}
\end{itemize}
As a general policy to simplify the notation we shall always drop
indices  when there in no danger of confusion.

We will generically assume that we have a closed, symmetric and semibounded from below form
$h$ with the dense domain $\q(h)\subset\H$ as given
in \cite[(VI.1.5)--(VI.1.11), pp. 308--310]{Kato76}.
The form $h$ which has a strictly positive lower bound will be called
\textit{positive-definite}. This is also a small departure from the terminology
of \cite[Section VI.2, pp. 310]{Kato76}.
Such form $h$ defines the self-adjoint and
positive definite operator $\mH$ in the sense of \cite[Theorem VI.2.23, pp. 331]{Kato76}. Furthermore, the operator
$\mH$ is densely defined with the domain $\d(\mH)\subset\q(h)$ and $\d(\mH^{1/2})=\q(h)$.
We also generically assume that $\mH$ has discrete eigenvalues
$
\lambda_1(\mH)\leq\cdots\leq\lambda_m(\mH)\leq\cdots<\les(\mH),
$
where we count the eigenvalues according to multiplicity. Another departure from the terminology of Kato is that
we use $h(\psi, \phi)$ to denote the value of $h$ on $\psi, \phi\in\q(h)$, but we write $h[\psi]:=h(\psi, \psi)$
for the associated \textit{quadratic form} $h[\cdot]$. We also
emphasize that we use $\cdot^*$ to denote the adjoint both in
the real as well as in the complex Hilbert space $\H$ as is customary in \cite[Chapters VI--VIII]{Kato76}.

Let us now fix our Rayleigh--Ritz terminology and outline the main
construction from \cite{Gru03_3,GruVes02}. We assume that we are
interested in approximating the eigenvalue $\lambda(\mH)$ of finite
multiplicity $m\in\N$.  Instead of only one test vector, as was the situation
in (\ref{TK-gruves}), we now need a \textit{test subspace} of
dimension $m$.

Let therefore $P$ be an orthogonal projection such that
$\di\ra(P)=m$ and $\ra(P)\subset\q(h)$. We call $\ra(P)$ the test
subspace for (the approximation of) $\lambda(\mH)$. The operator
$\Xi\in\lp(\ra(P))$, $\Xi=(\mH^{1/2}P|_{\ra(P)})^*\mH^{1/2}P|_{\ra(P)}$
will be called the (generalized) \textit{Rayleigh quotient}. Its
eigenvalues $\mu_1\leq\cdots\leq\mu_m$ will be called the
\textit{Ritz values} from the \textit{test subspace} $\ra(P)$ and
the vectors $u_i\in\ra(P)$, $\Xi u_i=\mu_i u_i$, $\|u_i\|=1$ will be called
the \textit{Ritz vectors}. We also define the operator
$\mW:\ra(P)^\perp\to\ra(P)^\perp$ as the one which is defined in
$\ra(P)^\perp$ by the form $h(P_\perp\cdot,P_\perp\cdot)$ in the
sense of \cite[Theorem VI-2.23, pp. 331]{Kato76}.

\subsection{A variant of the Temple--Kato inequality}
Before we outline the main form theoretic construction from
\cite{Gru03_3,GruVes02} let us illustrate our ``matrix theoretic''
approach to spectral theory by proving a variant of
(\ref{classical}) in a case when $\lambda_1(\mH)$ has a finite
multiplicity $m$.

Let $\mH$ be a self-adjoint operator which is bounded from below.
Let further, counting the eigenvalues according to multiplicity,
$\lambda:=\lambda_1(\mH)=\lambda_m(\mH)<\lambda_{m+1}(\mH)$,
$\lambda\ne0$ and we assume that we have a test subspace
$\ra(P)\subset\d(\mH)$, $\dim\ra(P)=m$. The environment space $\H$
can be decomposed as $\H=\ra(P)\oplus\ra(P)^\perp$ and $\mH$ can be
represented as a block-operator matrix\footnote{For more on block
operator matrices see for instance \cite{DavisKahWein82}.}
\begin{equation}
\mH=\begin{bmatrix}\Xi&K^*\\K&\mW\end{bmatrix},
\end{equation}
where $K\in\lp(\ra(P),\ra(P)^\perp)$, $K=P_\perp\mH P\big|_{\ra(P)}$.
Using the standard result\footnote{Also known as Kahan's residual
theorem in the case of the test subspace $\ra(P)$ of
Krylov-Weinstein inequality in the case of one test vector $\psi$,
see \cite{Chat83,Parlett80}.} \cite[Theorem 5.1]{DavisKahWein82} one
obtains that there exist $m$-eigenvalues
$\lambda_{j_1}(\mH)\leq\cdots\leq\lambda_{j_m}(\mH)$ such that
$$
|\lambda_{j_i}(\mH)-\mu_i|\leq\|K\|, \qquad i=1,\ldots m~.
$$
If we further assume that $\|K\|<\lambda_{m+1}(\mH)-\mu_m$, then
$$
|\lambda_{i}(\mH)-\mu_i|\leq\|K\|,\qquad i=1,\ldots m~.
$$
follows from the perturbation construction of \cite[Theorem 5.1 and
Remark 5.2]{DavisKahWein82}. Furthermore, the spectral calculus for
the self-adjoint operator  $\mW$ and the min-max formulae yield
$$
\|(\mW-\lambda\I)^{-1}\|\leq\frac{1}{\lambda_{m+1}(\mH)-\mu_m-\|K\|}.
$$
The assumption $\ra(P)\subset\d(\mH)$ allows us to justify the
following matrix representation
$$ \mH-\!\lambda\I=\!\!
\begin{bmatrix}\I&K^*\mB_{\textrm{abs}}^{-1}\\0&\I\end{bmatrix}\!\!
\begin{bmatrix}\Xi-\lambda\I-K^*\mB_{\textrm{abs}}^{-1}K&\\&\mB_{\textrm{abs}}\end{bmatrix}\!\!
\begin{bmatrix}\I&0\\\mB_{\textrm{abs}}^{-1}K&\I\end{bmatrix},\quad \mB_{\textrm{abs}}:=(\mW-\lambda\I).
$$
We can now use a generalization of the so called Wilkinson's trick
from \cite[pp. 183]{Parlett80}, which will be stated explicitly as
Theorem \ref{thm:wilkinson} below, to conclude that
$$
m=\dim\je(\mH-\lambda\I)=\dim\je(\diag(\Xi-\lambda\I-K^*\mB^{-1}_{\textrm{abs}}K,\mB_{\textrm{abs}}))
$$
and this can only happen if
\begin{equation}\label{eq:WilkonsonTrick}
\Xi-\lambda\I=K^*(\mW-\lambda\I)^{-1}K~.
\end{equation}
Furthermore, we establish
\begin{equation}\label{eq:TKSymm}
\tripleb\diag(\mu_i-\lambda)\tripleb\leq\frac{1}{\lambda_{m+1}(\mH)-\mu_m-\|K\|}
\tripleb K\tripleb\|K\|,
\end{equation}
where $\tripleb\cdot\tripleb$ is any unitary invariant norm. A
similar inequality holds for discrete eigenvalues which are located
in the interior of $\spec(\mH)$. In particular, if we denote the
(single) Ritz vector residuals by $r_i:=\mH u_i-\mu_i u_i$ and apply
the trace operator $\tr(\cdot)$ on (\ref{eq:WilkonsonTrick}), we obtain\footnote{
The inequality (\ref{eq:trace}) appeared with a better gap estimate in the bounded operator setting
in \cite{Ovtchinik1}, whereas the inequality (\ref{eq:TKSymm}) appeared in the matrix
setting in \cite{Sun}. The significance of this inequalities in this paper is to introduce the Schur complement technique
which will be the main tool later.}
\begin{equation}\label{eq:trace}
\sum^m_{i=1}|\mu_i-\lambda|\leq\frac{1}{\lambda_{m+1}(\mH)-\mu_m-\|K\|}\sum^m_{i=1}\|r_i\|^2.
\end{equation}
This generalizes the estimate from (\ref{TK:ineq}) to the case in
which $\lambda_1(\mH)$ has the multiplicity $m$. The quotient
$\frac{1}{\lambda_{m+1}(\mH)-\mu_m-\|K\|}$ is numerically inferior
to $\frac{1}{\lambda_{m+1}(\mH)-\mu_m}$, but (\ref{eq:TKSymm}) holds
for any unitary invariant norm. This is a mechanism which allows us to individually
treat Ritz vectors of different approximation properties. Furthermore, in view
of the discussion from \cite[pp. 8]{DavisKahan70} estimates
(\ref{eq:TKSymm}) give significant new information when compared
just with (\ref{eq:trace}). It could be argued that
$\frac{1}{\lambda_{m+1}(\mH)-\mu_m-\|K\|}$ is a reasonable
ingredient of the estimates since, in typical situations, one uses
\cite[Theorem 5.1 and Remark 5.2]{DavisKahWein82} to obtain a bound
for $\frac{1}{\lambda_{m+1}(\mH)-\mu_m}$. Also, see the discussion
on \cite[pp. 305]{Chat83}. For a way to compute $\tripleb K\tripleb$
see \cite[Section 9]{DavisKahan70}.

\subsection{The symmetric form approach}
The previous computation can not be justified in the case in which
$\ra(P)\subset\q(h)=\d(\mH^{1/2})$ but $\ra(P)\not\subset\d(\mH)$.
Precisely this is the case in which we are interested.

Let us now outline the main perturbation construction from
\cite{Gru03_3,GruVes02}. We start by defining the positive definite
form
$$
h_P(u, v)=h(Pu, Pv)+ h(P_\perp u, P_\perp v),\qquad u,v\in\q(h)
$$
and the self-adjoint operator $\mH_P$ which is defined by $h_P$ in
the sense of \cite[Theorem VI-2.23, pp. 331]{Kato76}. The operator
$\mH_P$ and the form $h_P$ are called the $P$-\textit{diagonal} part of $\mH$ and $h$, respectively.
We also define the form
$$
\delta h^P(u, v)=h(Pu, P_\perp v)+h(P_\perp u, P v), \qquad
u,v\in\q(h),
$$
which is an approximation defect in $\ra(P)$, since $\ra(P)$ is an
invariant subspace of $\mH$ if and only if $\delta h^P\equiv 0$, for
a proof see \cite{GruVes02}. Furthermore, it was shown in
\cite{Gru03_3,GruVes02} that
\begin{enumerate}
\item $\ra(P)$ reduces $\mH_P$
\item $\Xi=P\mH_PP\big|_{\ra(P)}$
\item $\ra(\mH^{-1}-\mH_P^{-1})$ is finite dimensional which implies that\\
$\spec_{ess}(\mH)=\spec_{ess}(\mH_P)$ according to the Weyl
theorem.
\end{enumerate}

The properties 1), 2) and 3) imply that $\mu_i\in\spec(\mH_P)$ and
 $\les(\mH)=\les(\mH_P)$ together with the assumption $\mu_i<\les(\mH)$
yields that $\mu_i$ are the eigenvalues of the operator $\mH_P$
with finite multiplicity. We are setting the scene for an
application of the relative perturbation theory from \cite[Chapters
VI--VIII]{Kato76} and so we will be able, regardless of the fact that\footnote{In
fact, it is even possible that $\d(\mH)\cap\d(\mH_P)=\{0\}$ and the
form approach is still applicable.} $\d(\mH)\ne\d(\mH_P)$, to interpret
$\mH$ as a perturbation of $\mH_P$ and thus bring $\mu_i$ in
connection with some component of $\spec(\mH)$.

This was the main line of argument in \cite{Gru03_3,GruVes02}.
Although some of the technical results about $\mH_P$, which we shall
now state are not explicitly given in \cite{Gru03_3,GruVes02} we
present them here without proof. However, all of their proofs are
obtainable as minor modifications of the arguments from
\cite{Gru03_3,GruVes02} and do not bring any new information.

Let us now look into the structure of this construction in more
detail. According to \cite[Theorem 4.5]{Gru03_3} the form
$
\delta h^P_s(\cdot, \cdot):=\delta h^P(\mH^{-1/2}_P\cdot,
\mH^{-1/2}_P\cdot)
$
defines the bounded operator $\delta H^P_s$ and
\begin{equation}\label{eq:deltaHs}
\|\delta H^P_s\|=\max_{\psi\in\ra(P)}\frac{(\psi,
\mH^{-1}\psi)-(\psi,\mH^{-1}_P\psi)}{(\psi,\mH^{-1}\psi)}.
\end{equation}
To examine $\delta H_s^P$ in further detail define
\begin{equation}\label{eq.sing_val_2}
\eta_i(P):=\Big[\max_{\substack{\mathcal{S}\subset\ra(P),\\
\dim(\mathcal{S})=m-i+1}}\min\big\{\frac{(\psi,
\mH^{-1}\psi)-(\psi,\mH^{-1}_P\psi)}{(\psi,\mH^{-1}\psi)}~
\big|~\psi\in\mathcal{S}, \|\psi\|=1\big\}\Big]^{1/2},
\end{equation}
for $i=1, \ldots, m$. Obviously, $\|\delta H^P_s\|=\eta_m(P)$. A
more detailed assessment of the prof of \cite[Theorem 4.1]{Gru03_3}
yields the following lemma.
\begin{lemma}\label{t:refinedGeom}
Let $h$ be positive definite, and let $\ra(P)$ be the test subspace
such that $\dim\ra(P)=m$. Assume further that
$$
\eta_{r+1}(P)\leq\eta_{r+2}(P)\leq\cdots\leq\eta_m(P)
$$
are all nonzero $\eta_i(P)$ from (\ref{eq.sing_val_2}). Then
$\eta_m(P)<1$ and  $\pm\eta_{r+1}(P), \ldots, \pm\eta_m(P)$ are all
non-zero eigenvalues of $\delta H_s^P$. Furthermore, $\eta_{r+1}(P),
\ldots, \eta_m(P)$ are all non-zero singular values of the operator  $K_s=
\big.\delta H_s^PP\big|_{\ra(P)}\in\lp(\ra(P),\ra(P)^\perp\!)$.
\end{lemma}
\begin{proof}
The proof of this lemma is implicitly contained in the proof of
\cite[Theorem 4.1]{Gru03_3}. We leave out most of the technical
details. We only explicitly present arguments that $\eta_m(P)<1$.
This fact was first established, in the matrix case, by Z. Drma\v{c},
Zagreb. Since $\ra(P)$ reduces $\mH_P$ we have
$$
h(\mH^{-1}_{P} f,v)=h_{P}(\mH^{-1}_{P} f,v)=(f,v),\qquad v\in\ra(P),
$$
i.e. $\mH^{-1}_{P}f$ is a Galerkin approximation from the subspace
$\ra(P)$ to $\psi=\mH^{-1}f$, which solves the problem $\mH\psi=f$.
With this in hand one computes
\begin{equation}\label{eq:galerkinB}
h(\mH^{-1} f-\mH^{-1}_{P} f,\mH^{-1} g-\mH^{-1}_{P}
g)=(f,\mH^{-1}g)-(f,\mH^{-1}_{P}g), \qquad f,g\in\ra(P).
\end{equation}
This implies
$(f,\mH^{-1}f)\geq(f,\mH_P^{-1}f)\geq 0$, $f\in\ra(P)$ and $0\leq\eta_i(P)\leq 1$ follows.
Assume $\eta_m(P)=1$, then there exists
$f\in\ra(P)\setminus\{0\}$ such that $(f,\mH^{-1}_Pf)=0$.
This is an obvious contradiction with the fact that $\mH_P$ is
positive definite.
\end{proof}

Now, set formally $\lambda_0(\mH):=0$, $\fG_0:=\infty$ and define
\begin{align}\fG_q&:=\min\Big\{
|\lambda_{q}(\mH)-\mu|\mu^{-1}:\mu\in\spec(\mH_P)\setminus\{\mu_1,...,\mu_m\}\Big\}\\
\gamma_s(\lambda_q)&:=\min\Big\{
\frac{\lambda_{q+m}(\mH)-\mu_m}{\lambda_{q+m}(\mH)+\mu_m},
\frac{\mu_1-\lambda_{q-1}(\mH)}{\mu_1+\lambda_{q-1}(\mH)}\Big\},\label{eq.uvijet123}
\end{align}
for $q\in\N$. This quantities---which measure the sensitivity of the eigenvalue $\lambda_q(\mH)$---will play a role in the
statement of the theorems in the next section. In the rest of the section we suppress the dependence of quantities
on $\mH$ and $P$ in the notation.

We now relate $\gamma_s(\lambda_q)$ and $\fG_q$ to $\eta_i $. The main result of
\cite{Gru03_3} established that given $\ra(P)\subset\q(h)$,
$\dim\ra(P)=m$ and $\mu_i<\les $ there exist $m$ eigenvalues
$\lambda_{j_1} \leq\cdots\leq\lambda_{j_m} $ such that
$$
\frac{|\lambda_{j_i} -\mu_i|}{\mu_i}\leq\eta_m ,\qquad i=1,
\ldots, m
$$
holds. Under an additional assumption on the location of the
unwanted component of the spectrum we can localize the approximated
eigenvalues and obtain that if e.g.
$\eta_m (1-\eta_m )^{-1}<(\lambda_{m+1} -\mu_m)(\lambda_{m+1} +\mu_m)^{-1}$
then
\begin{equation}\label{eq:broj}
\frac{|\lambda_{i} -\mu_i|}{\mu_i}\leq\eta_m ,\qquad i=1,
\ldots, m.
\end{equation}
This assumption is similar to the assumption of the Temple--Kato
inequality (\ref{TK:ineq}). For higher eigenvalues we have the
following variant of \cite[Theorems 5.1 and 5.2]{Gru03_3}, which we
present without proof.
\begin{lemma}\label{lemma:comb_1}
Let $\ra(P)$ be the test subspace for the positive definite form $h$
and let $\dim\ra(P)=m$ and $q\in\N$.
If $\frac{\eta_m}{1-\eta_m}<\gamma_s(\lambda_q)$
then $\fG_q>0$ and in particular\footnote{Here we assume
$\frac{c}{\infty}=0$, for $c>0$.}
\begin{align*}
\fG_q&\geq{\tiny \min\left\{
\frac{\mu_1(1-\eta_m )-(1+\frac{\eta_m }{1-\eta_m })\lambda_{q-1} }{
(1+\frac{\eta_m }{1-\eta_m })\lambda_{q-1} },
%&\qquad\qquad\qquad\qquad\qquad
\frac{(1-\frac{\eta_m }{1-\eta_m })\lambda_{q+m} -(1+\eta_m )\mu_m}{
(1-\frac{\eta_m }{1-\eta_m })\lambda_{q+m} }\right\}}.
\end{align*}
\end{lemma}

\section{Temple--Kato like inequality in the presence of Ritz value clusters}\label{Section3}

We now present the main contribution of this article. We will derive
relative eigenvalue estimates in the presence of Ritz value clusters.

In this section we will need to elaborate on the notion of the unitary invariant
operator norm (also known as symmetric or cross operator norms, cf.
\cite{Kato76,SimonTrace} and the references therein). This will
allow us to extract more information from $\ra(P)$ than what is
contained in $\eta_m(P)=\|K_s\|$. In this section we will be dealing
with only one orthogonal projection $P$, and so we simply write $\eta_i:=\eta_i(P)$, $i=1, \ldots, m$,
whenever there is no danger of confusion. Furthermore,
we write $\lambda_i:=\lambda_i(\mH)$, $i\in\N$ to simplify the notation.

 To say that the norm is \textit{unitary invariant} on $\SS\subset\lp(\H)$ means that, beside the
usual properties of any norm, it additionally satisfies:
\begin{description}
  \item[(i)] If $B\in\SS$, $A,C\in\lp(\H)$
  then $ABC\in\SS$ and
  \begin{equation}\label{eq:trple}
   \tripleb ABC \tripleb \leq\|A\| \tripleb B \tripleb \|C\|.
  \end{equation}
  \item[(ii)] If $A$ has rank $1$ then $ \tripleb A \tripleb =\|A\|$, where $\|\cdot\|$
  always denotes the standard operator norm on $\lp(\H)$.
  \item[(iii)] If $A\in\SS$ and $U,V$ are unitary on $\H$, then
  $UAV\in\SS$ and \begin{equation}\label{eq:unitary_inv_n}
  \tripleb UAV \tripleb = \tripleb A \tripleb.\end{equation}
  \item[(iv)] $\SS$ is complete under the norm $ \tripleb \cdot \tripleb $.
\end{description}
The subspace $\SS$ is defined as a $ \tripleb \cdot \tripleb
$--closure of the set of all degenerate operators in $\lp(\H)$. Such
$\SS$ is an ideal in the algebra $\lp(\H)$, cf. \cite{SimonTrace}.

A typical example of such a norm is the Hilbert--Schmidt norm
$\tripleb\cdot\tripleb_{HS}$. A bounded operator $H:\H\to\H$ is a
Hilbert--Schmidt operator if $H^*H$ is trace class and then, cf.
\cite[Ch. X.1.3]{Kato76},
$
\tripleb H \tripleb_{HS}:=
\sqrt{\tr H^*H}=[~\sum^\infty_{i=1}s_i(H)^2 ~]^{1/2}$.

Before we turn to the main theorem,
let us give an alternative---more common---definition for the approximation defects $\eta_i(P)$.
To this end we further exploit the Galerkin orthogonality
property of the $P$-diagonal part of $h$ and in particular the ramifications of
relation (\ref{eq:galerkinB}) from Lemma \ref{t:refinedGeom}. For any
$f\in\ra(P)$ we have
$$
h[\mH^{-1}
f-\mH^{-1}_Pf]=\Big[\sup_{\phi\in\q(h)\setminus\{0\}}\frac{|h(\mH^{-1}
f-\mH^{-1}_Pf,\phi)|}{h[\phi]^{1/2}}\Big]^{2}=\|\mH\Xi^{-1}
f-f\|_{\mH^{-1}}^2.
$$
With this
we can write (\ref{eq.sing_val_2}) in an
alternative form
\begin{equation}\label{eq.sing_val_3}
\eta_i(P):=\Big[\max_{\substack{\mathcal{S}\subset\ra(P),\\
\dim(\mathcal{S})=m-i+1}}\min\big\{\frac{\|\mH
\Xi^{-1}\psi-\psi\|_{\mH^ {-1}}^2}{\|\psi\|_{\mH^ {-1}}^2}~
\big|~\psi\in\mathcal{S}, \|\psi\|=1\big\}\Big]^{1/2},
\end{equation}
for $i=1, \ldots, m$.
\subsection{Operator matrices and the Wilkinson's trick}

As a first step we shall outline the Wilkinson's trick and state our
adaptation of this result as a theorem. This result yielded (\ref{eq:WilkonsonTrick}). We shall then
proceed to prove eigenvalue estimates. Let us now generalize the
Wilkinson's trick to operator matrices, cf. \cite[p. 183]{Parlett80}.

\begin{theorem}[Wilkinson's trick]\label{thm:wilkinson}
Let $A:\H_1\to\H_1$ and $X:\H_2\to\H_1$ be bounded operators and let
$A$ be self-adjoint. Assume further that $\mB:\H_2\to\H_2$ is
self-adjoint and that it has a bounded inverse and define
$
\mM=\begin{bmatrix} A&X\\X^*&\mB\end{bmatrix}
$,
to be understood as operator on $\H_1\oplus\H_2$. If ${\rm
dim}~\je(\mM)={\rm dim}~\H_1<\infty$ then
$$
A=X\mB^{-1}X^*.
$$
\end{theorem}
\begin{proof}
We shall adapt the Schur-complement technique from \cite[p.
183]{Parlett80}. Since $\mB^{-1}$ is assumed to be bounded we can
write
\begin{equation}\label{eq:cong}
\mM=\begin{bmatrix}\I&X\mB^{-1}\\0&\I\end{bmatrix}
\begin{bmatrix}
A-X\mB^{-1}X^*&0\\0&\mB\end{bmatrix}\begin{bmatrix}\I&0\\\mB^{-1}X^*&\I\end{bmatrix}=S
\mathbf{D} S^*.
\end{equation}
Both of the operator matrices
$$
S=\begin{bmatrix}\I&X\mB^{-1}\\0&\I\end{bmatrix},\qquad S^{-1}=\begin{bmatrix}\I&-X\mB^{-1}\\0&\I\end{bmatrix}
$$
define bounded operators on $\H_1\oplus\H_2$, and
so $\mathbf{D}=S^{-1}\mM S^{-*}$. This implies that $\d(\mM)=\d(\mathbf{D})$
and as a consequence of a simple dimension
counting we obtain that
\begin{align}
{\rm dim}~\je(\mM)&={\rm dim}~\je(\mathbf{D})<\infty\label{eq:wilinsonT2}
\end{align}
Since $\mB$ has a bounded inverse (\ref{eq:wilinsonT2}) can only be
true if $ A-X\mB^{-1}X^*=0$. This is the so called Wilkinson's trick
and it proves the statement of the theorem.
\end{proof}
\begin{remark}{\rm
Note that the theorem remains valid if we only assume that $\mB$ is
injective and $\mB^{-1}X^*$ is bounded. In this case we conclude that
$
A=X(\mB^{-1}X^*)
$. In the case when $\H_1$ is infinite dimensional the dimension counting cannot be used
to prove the result. Some spectral properties of Schur complements in a general situation can
be found in \cite{Tretter}.}\end{remark}

\begin{theorem}\label{prvo:t_DrmacHari}
Let $\ra(P)$ be the test subspace for the positive definite form $h$, as defined in Section \ref{section2},
and let $\dim\ra(P)=m$ and $q\in\N$. Assume that $\lambda_{q-1}<\lambda_q=\lambda_{q+m-1}<\lambda_{q+m}$
and $\frac{\eta_m}{1-\eta_m}<\gamma_s(\lambda_q)$ hold then
\begin{equation}\label{eq:realisticcase}
\tripleb\I-\lambda_q\Xi^{-1}\tripleb\leq
\frac{\eta_m}{\fG_q}\tripleb{\rm
diag}(\eta_1,\cdots,\eta_m)\tripleb~.
\end{equation}
In particular, for $\tripleb\cdot\tripleb=\tripleb\cdot\tripleb_{HS}$ and $\mu_i\in\spec(\Xi)$ we have the estimate
\begin{align}
%\label{eq:worstcase}
%\max_{i=1,...,m}\frac{|\lambda_{q}-\mu_i|}{\mu_i}&\leq\frac{1}{\fG_q}\eta_m^2\\
\Big[\sum^m_{i=1}\frac{(\lambda_q-\mu_i)^2}{\mu_i^2}\Big]^{1/2}&\leq\frac{\eta_m}{\fG_q}
\Big[\eta_1^2+\cdots+\eta_m^2\Big]^{1/2}.\label{eq:realisticcase2}
\end{align}
\end{theorem}
\begin{proof}
Let the form $h_P$ be the $P$-diagonal part of $h$. A modification of
\cite[Theorems 5.1 and 5.2]{Gru03_3} implies that
$$
h(\mH^{-1/2}_P\cdot,\mH^{-1/2}_P\cdot)-\lambda_q(\mH^{-1/2}_P\cdot,
\mH^{-1/2}_P\cdot)
$$
defines the bounded operator $H_s(\lambda_q)$, which allows the
operator matrix representation
\begin{equation}\label{we_see}
H_s(\lambda_q)=\left[\begin{matrix}\I-\lambda_q\Xi^{-1}&
K_s^*\\K_s&\I- \lambda_q\mW^{-1} \end{matrix}\right],
\end{equation}
with respect to $\H=\ra(P)\oplus\ra(P)^\perp$. Now, Lemma \ref{lemma:comb_1} implies that $\I-\lambda_q\mW ^{-1}$ is
invertible and we may use the Wilkinson's trick
%\footnote{The finite
%dimensional part of this proof is essentially contained in
%\cite{DrmHAri97}. Our analysis shows that the right
%setting for a generalization this result is the class of unbounded
%operators, rather then to attempt to work with the bounded (compact)
%operators. The reason is that in the infinite dimensional setting
%either the operator $\mH$ or the scaling $\mH^{-1}$ is
%going to be unbounded. We argue that proper choice is to chose the scaling as
%the ``well behaved'' object.}
to derive quadratic estimates (for some further technical details see Lemma
\ref{lemma:comb_1}). In particular we
have
$
\|(\I-\lambda_q\mW ^{-1})^{-1}\|=\frac{1}{\fG_q}<\infty
$.
Now temporarily set $\mB_{\textrm{rel}}:=(\I\!-\!\lambda_q\mW^{-1})$, then
\begin{align*}
H_s(\lambda_q)&=\!\!\left[\begin{matrix}\I &K_s^*\mB_{\textrm{rel}}^{-1}\\0&\I
\end{matrix}\right]\!\!\!
\left[\begin{matrix}(\I\!-\!\lambda_q\Xi^{-1})-K_s^*\mB_{\textrm{rel}}^{-1}K_s &0\\
0&\mB_{\textrm{rel}}\end{matrix}\right]
\left[\begin{matrix}\I & 0\\\mB_{\textrm{rel}}^{-1}K_s&\I\end{matrix}\right]
\end{align*}
and Theorem \ref{thm:wilkinson} yields
\begin{equation}\label{eq:apply}
\I-\lambda_q\Xi^{-1}=K_s^*(\I-\lambda_q\mW^{-1})^{-1}K_s.
\end{equation}
Property (\ref{eq:trple}) of a unitary invariant norm
$\tripleb\cdot\tripleb$ implies
$$
\tripleb \I-\lambda_q\Xi^{-1}\tripleb\leq\tripleb K_s\tripleb
\|(\I-\lambda_q\mW^{-1})^{-1}K_s\|.
$$
We apply (\ref{eq:unitary_inv_n}) and Theorem \ref{t:refinedGeom} on the last
inequality to complete the proof.
\end{proof}

The estimate of Theorem \ref{prvo:t_DrmacHari} was an equality up to
(\ref{eq:apply}). So, there is
more information  in (\ref{eq:apply}) than we have used so far.
\begin{theorem}\label{thm:lower}
If $\mu_m<\lambda_{m+1}$ and $\lambda_1=\lambda_m$ then
\begin{align}\label{eq:lowerL}
 \tripleb{\rm diag}(\eta_1^2, \cdots, \eta^2_m)\tripleb
&\leq\tripleb\I-\lambda_q\Xi^{-1}\tripleb\leq\frac{1}{\fG_1}
 \tripleb{\rm diag}(\eta_1^2, \cdots, \eta_m^2)\tripleb\\
\sum_{i=1}^{m}\eta_i^2&\leq\sum^m_{i=1}\frac{\mu_i-\lambda_i}{\mu_i}\leq
\frac{1}{\fG_1}
\sum_{i=1}^{m}\eta_i^2.\label{eq:TraceEst}
\end{align}
\end{theorem}
\begin{proof} The assumption $\mu_m<\lambda_{m+1}$ implies $\fG_1>0$. Now, this combined with
Lemma \ref{lemma:comb_1} yields \begin{equation}\label{eq:new}
K_s^*K_s\leq \I-\lambda_q\Xi^{-1}\leq\fG_1^{-1}K_s^*K_s.
\end{equation}
The conclusions now readily follow by an application of a norm $\tripleb\cdot\tripleb$ and the
trace operator $\tr(\cdot)$ on (\ref{eq:new}).\end{proof}

We now relate $\fG_1$ to the standard relative gap
which was extensively studied in the relative perturbation theory.
The result od Lemma \ref{lemma:comb_1} can be easily improved if we concentrate on the lower
part of the spectrum. This is a reasonable assumption.
\begin{corollary}\label{cor:reasonable}
Assume that
$\lambda_1 =\cdots=\lambda_{m} <\lambda_{m+1} $ and
$\eta_m <(\lambda_{m+1} -\mu_m)(\lambda_{m+1} +\mu_m)^{-1}$
then
$
\fG_1\geq\frac{\lambda_{m+1} -\mu_m}{\lambda_{m+1} +\mu_m}
$
and in particular
\begin{align}\label{eq:lowerL_C}
 \tripleb{\rm diag}(\eta_1^2, \cdots, \eta^2_m)\tripleb
&\leq\tripleb\I-\lambda_q\Xi^{-1}\tripleb\leq\frac{\lambda_{m+1} +\mu_m}{\lambda_{m+1} -\mu_m}
 \tripleb{\rm diag}(\eta_1^2, \cdots, \eta_m^2)\tripleb´.
\end{align}
\end{corollary}

\begin{remark}\label{rem:1}{\rm
Theorem \ref{thm:lower} illustrates why this ``relative'' or form
approach to the Temple--Kato inequality is more natural than the one
which yielded (\ref{eq:WilkonsonTrick}). The operator
$\mB_{\textrm{rel}}$ in the relative block operator
residual equation (\ref{eq:apply}) is such that $\I\leq\mB_{\textrm{rel}}^{-1}\leq\fG_1^{-1}\I$.
In comparison in the absolute block operator residual equation (\ref{eq:WilkonsonTrick}) it
only holds that $0\leq\mB_{\textrm{abs}}^{-1}$ and so no lower estimate is obtainable. In the more general case of
Theorem \ref{prvo:t_DrmacHari} the bounded operator $\mB_{\textrm{rel}}$
is indefinite but it is always boundedly invertible. For further discussion of the optimality of the
Schur complement approach see Section \ref{sec:exact}.}
\end{remark}

\begin{proposition}\label{prop:lower}
Let $\lambda_m<\lambda_{m+1}$ and let $\ra(P)$, $\dim\ra(P)=m$ be the test space for the form $h$ such that $2\eta_m<1$ then
$$
\frac{1\mu_1}{2\mu_m}\sum_{i=1}^m \frac{\|\mH u_i-\mu_i
u_i\|^2_{\mH^{-1}}}{\|\mH u_i\|^2_{\mH^{-1}}}\leq
\sum^m_{i=1}\frac{\mu_i-\lambda_i}{\mu_i}.
$$
\end{proposition}
\begin{proof}
Let $u_i\in\ra(P)$, $i=1,...,m$ be the Ritz vectors. The proof follows from
\begin{align*}
\sum^m_{i=1}\|\mH \frac{1}{\mu_i}u_i-u_i\|^2_{\mH^{-1}}&=\sum^m_{i=1}\Big[(u_i,\mH^{-1}u_i)-(u_i,\Xi^{-1}u_i)\Big]
=\sum^m_{i=1}\Big[(u_i,\mH^{-1}u_i)-\frac{1}{\mu_i}\Big]\\
&\leq\sum^m_{i=1}\Big[\frac{1}{\lambda_i}-\frac{1}{\mu_i}\Big]=\sum^m_{i=1}\frac{\mu_i-\lambda_i}{\lambda_i\mu_i}.
\end{align*}
 Note that according to (\ref{eq:Diag}) below, we have
$\sum_{i=1}^m\eta^2_i \leq \sum_{i=1}^m
\frac{\|\mH u_i-\mu_i u_i\|^2_{\mH^{-1}}}{\|\mH u_i\|^2_{\mH^{-1}}}$.
\end{proof}
\begin{remark}{\rm
The upper estimate in the setting of Proposition \ref{prop:lower}
can be achieved by a repeated application of the trace operator and
the estimate (\ref{eq:Diag}) to the identity (\ref{eq:apply}). The
estimate is rather technical an we leave it out. However, we
emphasize that we can recreate the framework of \cite[Proposition
2.3]{DrmacVeselic2} completely.}
\end{remark}

\subsection{A relationship with standard $\mH^{-1}$-norm residual estimates}
This section addresses the issue of the computability of $\eta_i(P)$ by relating these quantities to
the standard $\mH^{-1}$-norm estimates of the residuals associated to
the Ritz vector basis of $\ra(P)$. The proofs as well as the results
are technical and as such can be skipped on the first reading.
%Only the information from Lemma \ref{lemma:comb_1} below and
%the formula (\ref{eq.sing_val_3}) are necessary to follow the
%proof of the main results from Section \ref{Section3}.

Let us now reconsider the identities (\ref{eq.sing_val_2}) and (\ref{eq.sing_val_3}) and note that they
can be understood as generalized matrix eigenvalue problems. Assume
$u_1, \ldots, u_m$ are the Ritz vectors from $\ra(P)$, then for $i,j
=1, \ldots, m$, we define the matrices
\begin{align*}
\Omega_{ij}&=h(\mH^{-1}u_i-\Xi^{-1}u_i, \mH^{-1}u_j-\Xi^{-1}u_j)\\
\Psi_{ij}&=(u_i,\mH^{-1}u_j).
\end{align*}
Relation (\ref{eq:galerkinB}) from Lemma
\ref{t:refinedGeom} implies that $\Omega$ is a positive definite
matrix and in particular
\begin{align*}
\eta_i^2 &=\lambda_i(\Psi^{-1/2}\Omega\Psi^{-1/2})\\
\Omega_{ii}&=\|\mH\frac{1}{\mu_i}u_i-u_i\|^2_{\mH^{-1}}=\frac{1}{\mu_i^2}\|\mH
u_i-\mu_iu_i\|^2_{\mH^{-1}},\quad i=1, \ldots, m\\
D_\mu&\leq\Psi\leq(1+\fD_l)D_\mu,
\end{align*}
where $D_\mu=\diag(\mu_1^{-1}, \ldots, \mu_m^{-1})$ and $
\fD_l=\|D_\mu^{-1/2}(\Psi-D_\mu)D_\mu^{-1/2}\|$. Now, with the help
of $\mu_i=\|\mH u_i\|^2_{\mH^{-1}}$, we obtain
$$ \sum_{i=1}^m\lambda_i(D_\mu^{-1/2}\Omega
D_\mu^{-1/2})=\tr(D_\mu^{-1/2}\Omega D_\mu^{-1/2})=\sum_{i=1}^m
\frac{\|\mH u_i-\mu_i u_i\|^2_{\mH^{-1}}}{\|\mH u_i\|^2_{\mH^{-1}}},
$$
and so we conclude that
\begin{equation}\label{eq:preDiag}
\frac{1}{1+\fD_l}\sum_{i=1}^m
\Omega_{ii}\mu_i\leq\sum_{i=1}^m\eta^2_i \leq\sum_{i=1}^m
\Omega_{ii}\mu_i.
\end{equation}
Estimate (\ref{eq:preDiag}) can now be written as (cf. Proposition \ref{prop:lower})
\begin{equation}\label{eq:Diag}
\frac{1}{1+\fD_l}\sum_{i=1}^m \frac{\|\mH u_i-\mu_i
u_i\|^2_{\mH^{-1}}}{\|\mH
u_i\|^2_{\mH^{-1}}}\leq\sum_{i=1}^m\eta^2_i \leq \sum_{i=1}^m
\frac{\|\mH u_i-\mu_i u_i\|^2_{\mH^{-1}}}{\|\mH u_i\|^2_{\mH^{-1}}}.
\end{equation}
By a
similar argument one can conclude that asymptotically (as $P$
converges to the $m$-dimensional spectral subspace) we have as a heuristic
$$
\frac{\|\mH u_i-\mu_i u_i\|_{\mH^{-1}}^2}{\|\mH
u_i\|^2_{\mH^{-1}}}\sim\eta_{m-i+1}^2(P).
$$
This indicates that $\eta_i(P)$ represent a canonical choice of
residuals from $\ra(P)$, i.e. not defined by the Ritz vectors $u_i$
but rather the vectors which are selected by the variational formulae
(\ref{eq.sing_val_3}).

\begin{remark}\label{rem:Golub}{\rm
The definition of $\eta_i $ indicates that the problem of
computing (or estimating) $\eta_i $ requires the solution of the
$m\times m$ positive definite generalized eigenvalue problem. Since
$m$ is the multiplicity of the eigenvalue of interest, the
computational cost of the solution of such problem is negligible.
The main problem is how to evaluate or estimate the moments $(u_i,
\mH^{-1}u_j)$, $i,j=1, \ldots, m$ without actually inverting the
operator $\mH^{-1}$. For some possibilities to do this see
\cite[Section 3.]{Arioli}, \cite[Section 5.]{Golub} or \cite[Remark
8]{RitzDrm96}.}
\end{remark}
\section{A simple non-inhibited stiff problem}\label{Section5}
These estimates have been used in \cite{GruPhd} to study a
 class of eigenvalue problems which is given by the family of
positive definite forms
\begin{equation}\label{e_prva}
h_\kappa(u,v)=h_b(u,v)+\kappa^2h_e(u,v),\quad \kappa~\text{ large }.
\end{equation}
The forms $h_b$ and $h_e$ are assumed to be symmetric, closed and nonnegative and we further assume that
$h_b+h_e$ is positive definite in $\H$ and that $\q(h_b+h_e)$ is dense in $\H$.  Family (\ref{e_prva}) can
always be considered as a perturbation of $h_b+h_e$ (after an
obvious change of variable $\kappa$) rather than as a perturbation
of $h_b$ and so we assume, without affecting the generality of
results, that $h_b$ is positive definite and densely defined.

A detailed study of the spectral property of the families like
(\ref{e_prva}) is beyond the scope of this article and will be
reported in subsequent publication. We will now consider a very simple
problem of this form, and note that (\ref{e_prva}) motivated the
example (\ref{OvtEx}). Let $H^1_0[0,1]$ and $H^1_0(\R_{+})$, $\R_{+}:=\big[0,\infty\big>$ be the
standard Sobolev spaces. We also identify the functions from
$H^1_0[0,1]$ with their extension by zero to the whole of $\R_{+}$
and write $ H_0^1[0,1]\subset H_0^1(\R_{+})$.
Consider the family of positive definite
forms
\begin{equation}\label{e_schrod}
h_\kappa(u,v)=\int_0^\infty u'v'~dx + \kappa^2\int_1^\infty uv~dx,
\quad u,v\in H^1_0(\R_+).
\end{equation}
By $\mH_\kappa$ we denote the positive definite operator defined by
$h_\kappa$ in (\ref{e_schrod}). The operators $\mH_\kappa$ converge
in the generalized sense to the operator $\mH_\infty$, which is
defined by the form $h_\infty(u,v)=\int^1_0 u'v'~dx$, $u,v\in
H^1_0\left[0,1\right]$. For further details on this convergence see \cite{GruPhd} and the references therein.
We also formally write
$\mH_\kappa=-\partial_{xx}+\kappa^2\chi_{\left[1,\infty\right>}$ and
$\mH_\infty=-\partial_{xx}$.  As a test function(s)
we chose
\begin{equation}\label{drugo:e_testfunkcija}
u_q(x)=\begin{cases}\sqrt{2}\sin(k \pi x),&0\leq x\leq 1\\
0,& 1\leq x\end{cases}~, q\in\N .
\end{equation}
Note that here $u_q\in\q(h_\kappa)$ but $u_q\not\in\d(\mH_\kappa)$.
The eigenvalues of the operator $\mH_\kappa$ have to be described
implicitly. Let $ \mH_\kappa v^\kappa=\lambda^\kappa v^\kappa, $
then $v^\kappa\in C^1(\R_+)$ is
$$
v^\kappa(x)=\begin{cases}\sin(\sqrt{\lambda^\kappa}x),&0\leq x\leq 1\\
\frac{\sin\sqrt{\lambda^\kappa}}{e^{-\sqrt{\kappa^2-\lambda^\kappa}}}~e^{-\sqrt{\kappa^2-\lambda^\kappa}~x},&1\leq
x\end{cases}
$$
and $\lambda^\kappa$ is a solution of the equation
\begin{equation}\label{trece:nonlinear}
\sqrt{\kappa^2-\lambda^\kappa}=-\sqrt{\lambda^\kappa}\cot(\sqrt{\lambda^\kappa}).
\end{equation}
The quotient
$\frac{\lambda^\infty_1-\lambda_1^\kappa}{\lambda^\infty_1}$ can be
represented (for $\kappa\to\infty$) by a convergent Taylor series
\begin{equation}\label{drugo:e_razvoj}
\frac{\lambda^\infty_1-\lambda_1^\kappa}{\lambda^\infty_1}=
2\frac{1}{\kappa}-3\frac{1}{\kappa^2}+8\left(\frac{1}{2!}+\frac{1}{4!}\pi^2\right)
\frac{1}{\kappa^3}-10\left(\frac{1}{2!}+\frac{4}{4!}\pi^2\right)\frac{1}{\kappa^4}+\cdots~.
\end{equation}
We directly compute $\eta^2_\kappa(u_q):=\frac{2}{3+\kappa}$
%\footnote{By $(\cdot, \cdot)_{L^2(\R_{+})}$ and $(\cdot,
%    \cdot)_{L^2[0, 1]}$ we denote the scalar products in $L^2(\R_{+})$ and
%$L^2[0, 1]$ respectively.}
%$$
%\eta^2(u_q)= \frac{(u_q, \mH_\kappa^{-1} u_q)_{L^2(\R_{+})}-(u_q,
%\mH_\infty^{-1} u_q)_{L^2[0,1]}}{(u_q, \mH_\kappa^{-1}
%u_q)_{L^2(\R_{+})}}=\frac{2}{3+\kappa}
%$$
and combine it with (\ref{trece:nonlinear}) and the first order
estimate from (\ref{eq:broj}) to obtain
$$
\big(1-\sqrt{\frac{2}{3+\kappa}}\big)4\pi^2=:D(\kappa)\leq\lambda_2(\mH),\qquad
\kappa\geq 5.
$$
Theorems \ref{prvo:t_DrmacHari} and \ref{thm:lower} now yield
\begin{equation}\label{eq:applying}
\frac{2}{3+\kappa}\leq\frac{\lambda^\infty_1-\lambda_1^\kappa}{\lambda^\infty_1}\leq\frac{D(\kappa)+\pi^2}{D(\kappa)-\pi^2}\frac{2}{3+\kappa}=
\frac{10}{3\kappa}+\frac{1}{\sqrt{\kappa}}O\Big(\frac{1}{\kappa}\Big),
\qquad \kappa \geq 5,
\end{equation}
which is a tight estimate on the behavior of
$\frac{\lambda^\infty_1-\lambda_1^\kappa}{\lambda^\infty_1}$.
Similar estimates hold for other eigenvalues and eigenvectors, too.
This example illustrates the ``efficiency'' of this \textit{a
posteriori} estimator. Furthermore, it indicates a role which is
played by the first order estimates from \cite{Gru03_3} in the
general theory. For some details of the computation see
\cite{GruPhd}. The Schroedinger operators in higher dimensions have
also been studied in \cite{GruPhd}. The estimate for $\eta^2(u_q)$
can in this case be computed by a use of the advanced probabilistic
techniques from \cite{DemuthJeskeKirsch} or by a use of the boundary
layer techniques from \cite{BruneauCarbou} (naturally, under the
assumption that the domain is finite).

\subsection{A framework for proving the asymptotic exactness}\label{sec:exact}
Let us go back to Remark \ref{rem:1}. The conclusion of Proposition
\ref{prop:lower} does not appear to be completely satisfactory. The
factor $\mu_1/\mu_m$ limits its applicability to a couple of the
lowermost eigenvalues of $\mH$. The true power of the Schur
complement technique can be seen if we rewrite (\ref{eq:apply}) as
\begin{equation}\label{eq:exact}
\I-\lambda_q\Xi^{-1}=K_s^*K_s+\lambda_qK_s^*\mW^{-1/2}(\I-\lambda_q\mW^{-1})^{-1}\mW^{-1/2}K_s.
\end{equation}
After applying the trace operator on (\ref{eq:exact}) and utilizing Lemma \ref{t:refinedGeom} we obtain
$$
\frac{\sum^m_{i=1}\frac{\mu_i-\lambda_q}{\mu_i}}{\sum^m_{i=1}\eta^2_i}=1+\frac{\tr(
\lambda_qK_s^*\mW^{-1/2}(\I-\lambda_q\mW^{-1})^{-1}\mW^{-1/2}K_s)}{\sum^m_{i=1}\eta^2_i}.
$$
In our $\kappa$ dependent problem we use this to prove (cf. (\ref{drugo:e_razvoj}) and (\ref{eq:applying}))
$$\lim_{\kappa\to\infty}
\frac{\frac{\lambda^\infty_q-\lambda_q^\kappa}{\lambda^\infty_q}}{\eta^2_\kappa(u_q)}=1.$$
Furthermore, we see why this convergence is pretty rapid. In a general situation we perform this analysis by
comparing the
singular values $s_i(K_s(\kappa))$ with $s_i(\mW^{-1}_\kappa K_s(\kappa))$ and
noticing that $s_i(\mW^{-1}_\kappa K_s(\kappa))$
is of higher order in $\kappa^{-1}$. Here we have assumed an obvious modification of the block matrix
representation (\ref{we_see}) for the $\kappa$ dependent problem.
Exploiting (\ref{eq:exact}) in the general setting of (\ref{e_prva})
as well in as in the setting of finite element approximations is
beyond the scope of this paper and is a subject of subsequent
reports.

\section{Finite element computations}\label{Section6}
As a further explicitly solvable model example let us consider the
family of eigenvalue problems
\begin{eqnarray}
\nonumber -\psi''-\alpha \,\psi&=&\omega~ \psi,\\
e^{{\rm i} \theta}\psi(0)&=&\psi(2\pi),\label{prvo2:e_problem1}\\
\nonumber e^{{\rm i} \theta}\psi'(0)&=&\psi'(2\pi),
\end{eqnarray}
where $\theta\in\left[0, \pi\right]$ and we chose $\alpha\in\R$ so
that the eigenvalues remain positive. The weak formulation of
(\ref{prvo2:e_problem1}) is given by
$$
h(\psi, v_i)=\lambda_i(\mH)(\psi,v_i),\qquad \psi\in\q(h),
$$
where
\begin{align}\label{prvo:e_formazaH}
h(\psi, \phi)&:=\int^{2\pi}_0\big(\overline{\psi'}\phi'-\alpha \overline{\psi}\phi\big)~,
&\q(h):= \{\psi~|~\psi,\psi'\in \H, e^{{\rm
i}\theta}\psi(0)=\psi(2\pi)\}
\end{align}
and $\H=L^2[0,2\pi]$.
The eigenvalues of the problem (\ref{prvo2:e_problem1}) as well as
the Green function of the operator $\mH$, which is defined by
(\ref{prvo:e_formazaH}) are explicitly known, see \cite[Theorem
XIII.89, Volume 4. pp. 293]{ReedSimonSvi} and \cite[Equation
(XIII.154), pp. 292]{ReedSimonSvi}. In particular we have
\begin{align}
\lambda_1(\mH)&=\big(-1+\frac{\theta}{2\pi}\big)^2 -\alpha,
\;\;\lambda_2(\mH)=\big(\frac{\theta}{2\pi}\big)^2 -\alpha,\;\;
\lambda_3(\mH)=\big(1+\frac{\theta}{2\pi}\big)^2 -\alpha,\\
v_{1}(t)&=e^{-\imag \big(-1+\frac{\theta}{2\pi}\big)t},
\phantom{-\alpha,\;\lambda\;,,} v_{2}(t)=e^{-\imag \big(\frac{\theta}{2\pi}\big)t},\;
\phantom{-\alpha,}v_{3}(t)=e^{-\imag \big(1+\frac{\theta}{2\pi}\big)t}
\end{align}
and\footnote{We implicitly assume that $\H=L^2[0, 2\pi]$.}
\begin{align*}
(\psi, \mH^{-1}\phi)&=\int^{2\pi}_0 \,d t_1\int^{2\pi}_0 G(t_1-
t_2)\overline{\psi(t_2)}\phi(t_1)\,d t_2\\
G(t_1-t_2)&= \frac{\imag }{2\sqrt{\alpha }}\,\left( e^{\imag
\,{\sqrt{\alpha }}\,|t_1 - t_2|} +
      \frac{e^{\imag \,\left( t_1 - t_2 \right) \,{\sqrt{\alpha }}}}
       {-1 + e^{-2\,\imag \,\pi \,{\sqrt{\alpha }} - \imag \,\theta }} +
      \frac{e^{\imag \,\left( t_2-t_1 \right) \,{\sqrt{\alpha }}}}{-1 + e^{-2\,\imag \,\pi \,{\sqrt{\alpha }} + \imag \,\theta }}
      \right).
\end{align*}
Let us now choose $\theta=\pi$ and $\alpha=0.2499$ for our numerical
experiment. With this choice of parameter the problem
(\ref{prvo2:e_problem1}) is almost singular  and
$\lambda_1(\mH)=\lambda_2(\mH)$. For $N\in\N$ define the finite
element space
\begin{align*}
\vp^1_{N}&=\left\{\psi~|~\psi\in C[0,2\pi],-\psi(0)=\psi(2\pi),
%\right.\\
%&\left.\qquad\qquad
\psi\hbox{ is linear in } \mathcal{I}_p, p=1,\ldots, N\right\},
\end{align*}
where $\mathcal{I}_p:=\left<\frac{(p-1)
2\pi}{N},\frac{p~ 2\pi}{N}\right>$, and use
\begin{equation}\label{eq:RR}
\mu_i(\vp^1_N):=\max_{\substack{\mathcal{S}\subset\vp^1_{N}\\ {\rm
dim}\mathcal{S}={\rm dim}\vp^1_{N} -
i}}\min_{\psi\in\mathcal{S}\setminus\{0\}}\frac{h(\psi,\psi)}
{(\psi,\psi)},\qquad i=1,2
\end{equation}
to define the Rayleigh-Ritz approximations to the eigenvalue
$\lambda_1(\mH)=\lambda_2(\mH)$. Let also $u_i(\vp^1_N)\in\vp^1_N$,
$i=1,2$ be two vectors of norm one for which
$\mu_i(\vp^1_N)=h[u_i(\vp^1_N)]$, $i=1,2$ holds. Now, let
$P(\vp^1_N)$ be an orthogonal projection onto the linear span of
$\{u_1(\vp^1_N), u_2(\vp^1_N)\}$ and set
$\Xi_{P(\vp^1_N)}=\mH_{P(\vp^1_N)}P(\vp^1_N)$. We now apply Theorems
\ref{prvo:t_DrmacHari} and \ref{thm:lower} on the projections
$P(\vp^1_N)$ and display the results on Table \ref{f_femmod2}.

\begin{table}
{\small
\begin{center}
\begin{tabular}{|c|c|c|c|}\hline
  \begin{minipage}[c]{1cm}\begin{center}$~$\\[0.5em]N\\[1em]\end{center}
  \end{minipage}&\begin{minipage}[c]{2.5cm}\begin{center}
   $$\text{estimate (\ref{eq:lowerL})}$$
  \end{center}\end{minipage}&\begin{minipage}[c]{2.5cm}\begin{center}
  $$\tripleb \I-\lambda\Xi_{P(\vp^1_N)}^{-1}\!\tripleb_{HS}$$
  \end{center}\end{minipage}&
  \begin{minipage}[c]{2.2cm}\begin{center}
  $$\text{estimate (\ref{eq:realisticcase2})}$$
  \end{center}\end{minipage}\\[2em]\hline
  &&&\\
  40  & 7.9540e-001& 7.9540e-001 & 7.9558e-001  \\
  60  & 5.1413e-001& 5.1413e-001 & 5.1422e-001  \\
  80  & 3.4389e-001& 3.4389e-001 & 3.4393e-001  \\
  100 & 2.4120e-001& 2.4120e-001 & 2.4123e-001  \\
  120 & 1.7671e-001& 1.7671e-001 & 1.7673e-001    \\
 &&&\\\hline
\end{tabular}
~\\[1em]
\caption{The performance of the estimates (\ref{eq:lowerL}) and
(\ref{eq:realisticcase2}) on the family of test spaces
$\ra(P(\vp^1_N))$ and for the choice of the norm
$\tripleb\cdot\tripleb=\tripleb\cdot\tripleb_{HS}$. The
computational details can be found in \cite[Section 2.7.3, pp.
64]{GruPhd}.}\label{f_femmod2}
\end{center}}
\end{table}
\subsection{Hierarchical error estimation}

The results from Table \ref{f_femmod2} show that
$\eta_i,\ldots,\eta_m$ accurately capture the behavior of the relative
error as $\frac{1}{N}\to 0$. The explicit knowledge of the Green
function is most certainly an information which cannot in general be
assumed when considering higher dimensional eigenvalue problems. Let
us now consider an application of these estimates in the context of
the adaptive finite element methods for divergence type elliptic
self-adjoint operators in dimension two. We only present a feasibility
argument, an algorithmic development will be a subject of a subsequent report.

For the sake of definiteness let $\H=L^2(\region)$, where
$\mathcal{R}$ is assumed to be a bounded \textit{polygonal domain}
and let
\begin{equation}\label{eq:form}
h(u,v)=\int_{\region}(\nabla u)^*~\nabla v, \qquad u,v\in
\q(h)=H^1_0(\region).
\end{equation}
By $H^1(\mathcal{R})=\{u\in \H:\norm{u}^2+\|\nabla
u\|^2<\infty\}$ we denote the standard first order \textit{Sobolev
space}. The gradient $\nabla$ is meant in the weak sense and
$\|\cdot\|$ denotes the norm on $L^2(\mathcal{R})$ and
$H^1_0(\mathcal{R})\subset H^1(\mathcal{R})$ is assumed to be
equipped with the norm $\norm{u}_{E}:=\|\nabla u\|=h[u]^{1/2}$ and
it consists of those $H^1(\mathcal{R})$ functions which vanish on
the boundary of $\mathcal{R}$ in the sense of the trace operator.

The set $\mathcal{T}_d$ is called a \textit{triangulation of the
polygonal domain} $\region$ if it consists of the triangles such
that union of these triangles is $\overline{\mathcal{R}}$ and such
that the intersection of two such triangles either consists of a
common side or of a common vertex of both triangles or is empty. By
$d$ we denote the maximal diameter of all triangles in
$\mathcal{T}_d$. For a given triangulation $\mathcal{T}_d$ we define
the finite dimensional function spaces:
\begin{align*}
\vp^1_{d}&=\{u\in \q~|~ u\in~C(~\overline{\mathcal{R}}~) \text{ and } \left.v\right|_K \textrm{ is
linear }, K\in\cT_d\}\\
\vp^2_{d}&=\{u\in \q~|~u\in~C(~\overline{\mathcal{R}}~) \text{ and } \left.v\right|_K \textrm{ is
quadratic }, K\in\cT_d\}
\end{align*}
and the orthogonal projections $V_{i,d}$, $i=1,2$ such that
$\ra(V_{i,d})=\vp^i_{d}$, $i=1,2$. To simplify the notation we write
$\mH_{i,d}:=\mH_{V_{i,d}}$, $i= 1, 2$ and also define the orthogonal projection
$P_d$ such that $\ra(P_d)$ equals the linear span of $\{u_{1,d}, u_{2,d}\}$. In what follows we assume, as in
\cite{Nochetto2002}, that $\mathcal{T}_d$ is
graded and shape regular family of triangulations and that it satisfies the nondegeneracy
property \cite[Assumption 4.1]{Nochetto2002}.  Let us assume that we
have
\begin{equation}\label{eq:RR2}
\mu_{i, d}:=\max\big\{\min_{\psi\in\mathcal{S}\setminus\{0\}}\frac{h(\psi,\psi)}
{(\psi,\psi)}:\mathcal{S}\subset\vp^1_d, {\rm dim}\mathcal{S}={\rm
dim}\vp^1_d - (i+1)\big\},\quad i=1,2
\end{equation}
and $u_{i, d}\in\vp^1_{\cT_d}$, $i=1,2$ are chosen so
that\footnote{This can be checked by a direct computation.}
$\mH_{1,d}~u_{i, d}=\mu_{i, d}~u_{i, d}$, $i=1,2$.
The result \cite[Theorem 1.1]{Nochetto2002} and in particular the
last remark on \cite[pp. 12]{Nochetto2002} yield the estimate
\begin{equation}\label{eq:saturation}
\frac{h[\mH^{-1}u_{i, d}-\mH_{2,d}^{-1}u_{i, d}]^{1/2}}
{h[\mH^{-1}u_{i, d}-\mH_{1,d}^{-1}u_{i, d}]^{1/2}}=
\frac{\|\mH^{-1}u_{i, d}-\mH_{2,d}^{-1}u_{i, d}\|_E}
{\frac{1}{\mu_{i, d}}\|\mH u_{i, d}-
\mu_{i, d}u_{i, d}\|_{\mH ^{-1}}}\leq\alpha,\qquad
i=1,2
\end{equation}
with the constant $\alpha$ which depends solely on the shape
regularity of $\cT_d$. Set
$r_{i, d}=\mH_{2,d} u_{i, d}-\mu_{i, d}
u_{i, d}$. Combining (\ref{eq:saturation}) and
\cite[Estimate (2.16)]{MR1231320} we conclude that there exists constants $C_*$ and $c_*$, solely depending on the shape regularity of $\cT_d$, such that
\begin{equation}\label{eq:trading}
c_*\frac{\|r_{i, d}\|^2_{\mH_{2,d}^{-1}}}{\mu_{i, d}}\leq\frac{\|\mH
u_{i, d}-\mu_{i, d}
u_{i, d}\|^2_{\mH^{-1}}}{\|\mH
u_{i, d}\|^2_{\mH^{-1}}}\leq
C_*\frac{\|r_{i, d}\|^2_{\mH_{2,d}^{-1}}}{\mu_{i, d}}, \qquad i=1, 2.
\end{equation}
This estimate can now be directly
plugged into the trace type estimates from Theorem \ref{thm:lower} or
Proposition \ref{prop:lower}. Furthermore, Remark \ref{rem:Golub} allows us to exploit
other unitary invariant norms with similar ease.

\begin{remark}{\rm
Note that $\mH_{2,d}|_{\vp^2_{d}}\in\mathcal{L}(\vp^2_{d})$.
The $\mH_{2,d}^{-1}$ norms of the residuals
$r_{i, d}\in\vp^2_{\cT_d}$ can efficiently (cheaper than
when solving a linear system) be approximated as functions of the
vectors $w_{i, d }\in\vp^2_{\cT_d}\ominus\vp^1_{\cT_d}$,
$i=1,2$, as given by \cite[Theorem 2.1]{MR1231320}. Similar
consideration has been explored in \cite[Estimates (29)--(30)]{Neymeyr02}
(cf. \cite[Theorem 2.2]{MR1231320}), but in comparison our estimates
give more explicit information on the dependence of the constants on
the mesh and provide the optimality argument, too.}
\end{remark}

To summarize, the arguments of Remark \ref{rem:Golub} indicate that it is
possible to estimate the $\mH^{-1}$ norm of the residual cheaper than it takes to solve the linear
system. Furthermore, we have shown that when deciding on the convergence of
the finite element method the size of
$\|r_{i, d}\|^2_{\mH_{2,d}^{-1}}/\mu_{i, d}$
should be compared to the relative gap measure $\fG_q$ to decide if the
approximation is good enough. By this we mean if the whole multiplicity of the
target eigenvalue has been resolved by $\ra(P_d)$. To get a feeling for this statement one should remember the
picture of \cite[Example 2.1]{DrmacVeselic2}.

\section{Conclusion}
The main benefit of our approach is that, as the theoretical
considerations from Section \ref{Section5} and Table \ref{f_femmod2}
corroborate, up to (\ref{eq:trading}) we have had globally optimal
estimates for the eigenvalue error (i.e. almost no information was
lost). After (\ref{eq:trading}) we have started aggressively trading
off accuracy for speed. Our theory is such that this
can be achieved, in numerous situations, by a simple combination of
the Galerkin orthogonality condition and any of the ``of the shelf''
results like
the those from \cite{MR1231320,Nochetto2002} or \cite{BruneauCarbou,DemuthJeskeKirsch} in the
singularly perturbed setting. On top of this comes
the heuristic insight from \cite[Section 3]{DrmacVeselic2} which
indicates that we have properly identified the components of the
error as given by (\ref{general_form}). It should be noted that our
Theorem \ref{prvo:t_DrmacHari} directly corresponds to
\cite[Proposition 2.3]{DrmacVeselic2}, since both are motivated by
\cite{DrmHAri97}. Furthermore, the estimates from Theorems
\ref{prvo:t_DrmacHari} and \ref{thm:lower} can be combined with
\cite[Theorem 6.1]{Gru03_3} and the well known identity from
\cite[Ad (v), pp.617]{Larson2000} to obtain optimal estimates for
the eigenvalue error $\|u_i-v_i\|_E/\|v_i\|_E$. The estimates
can even be obtained in a situation in which the multiple eigenvalue
splits in a cluster of eigenvalues. This is a subject of the
followup report. As a conclusion let us remember the remarks 1), 2), 3) from the
Introduction. We have introduced a matrix analytic techniques which
tackle both test vectors outside the domain of definition of the
operator and the multiplicity of the approximated eigenvalue in a
natural and constructive way. Furthermore, Remark \ref{rem:Golub} opens a way to exploiting other
unitary invariant norms for scaling robust eigenvalue estimation.
\section*{Acknowledgement}The author would like to
thank Prof. Dr. Kre\v{s}imir Veseli\'{c}, Hagen and Prof. Dr. Volker Enss,
Aachen for helpful discussions and support during the research and the
preparation of this manuscript. The author also thanks Dr. Mario Arioli, Didcot
for a helpful discussion and for pointing out the reference \cite{Golub}.

\bibliographystyle{abbrv}
%\bibliography{../bibliografija}
%\end{document}
\newpage
\def\cprime{$'$} \def\cprime{$'$} \def\cprime{$'$}

\end{document}